\documentclass[authoryear, 3p]{elsarticle}

 \usepackage{amssymb}
 \usepackage{amsthm}
 \usepackage{amsmath}
 \usepackage{natbib}
 \usepackage{tikz}
 \usepackage{graphicx}
 \usepackage{epstopdf}
 \usepackage{xurl}
 \usepackage{hyperref}
 \usepackage{multirow}
 \usepackage{lineno}
 \theoremstyle{plain}
 \newdefinition{definition}{Definition}
\newtheorem{theorem}{Theorem}

  \newtheorem{proposition}{Proposition}
 \newdefinition{remark}{Remark}
 \newdefinition{property}{Property}
  \newdefinition{example}{Example}

\newcommand{\pp}{\mathbb{P}}
\newcommand{\xx}{\tilde{X}}
\newcommand{\yy}{\tilde{Y}}
\renewcommand{\ss}{\tilde{s}}
\renewcommand{\tt}{\tilde{t}}
\DeclareMathOperator{\sign}{sign}

\journal{arXiv}

\begin{document}

\begin{frontmatter}



\title{Kendall's tau estimator for bivariate zero-inflated count data}


\author[TUe]{Elisa PERRONE\footnote{Corresponding author, email \url{e.perrone@tue.nl}}}
\author[TUe]{Edwin R. van den HEUVEL}
\author[TUe]{Zhuozhao ZHAN}

\address[TUe]{Department of Mathematics and Computer Science, Eindhoven University of Technology, Groene Loper 5, 5612 AZ, Eindhoven, The Netherlands}

\begin{abstract}
In this paper, we extend the work of \cite{pimentel_association_2015} and propose an adjusted estimator of Kendall's $\tau$ for bivariate zero-inflated count data. We provide achievable lower and upper bounds of our proposed estimator and show its relationship with current literature.
In addition, we also suggest an estimator of the achievable bounds, thereby helping practitioners interpret the results while working with real data. The performance of the proposed estimator for Kendall's $\tau$ is unbiased with smaller mean squared errors compared to the unadjusted estimator of \cite{pimentel_association_2015}. Our results also show that the bound estimator can be used when knowledge of the marginal distributions is lacking.
\end{abstract}

\begin{keyword}
Kendall's tau \sep Bivariate zero-inflated count data \sep Fr\'echet-Hoeffding bounds
\end{keyword}

\end{frontmatter}

\section{Introduction}
Zero-inflated data naturally appears in many applications such as health care and ecology \citep{Moulton_1995,Arab_2012}. 
Analyzing zero-inflated data is challenging as the high amount of observations in zero invalidates standard statistical techniques. 
For example, assessing the level of dependence between two zero-inflated random variables becomes a difficult task as standard rank-based measures of association such as Kendall's $\tau$ and Spearman's $\rho$ cannot be applied directly due to the large amount of tied values in zero making any tie-breaking adjustment unsatisfactory \citep{hollander_13}.
The importance of deriving accurate estimators of popular association measures, such as Kendall's $\tau$ and Spearman's $\rho$, for bivariate zero-inflated distributions motivated recent work on the topic. In \cite{pimentel_kendalls_nodate} and \cite{pimentel_association_2015}, the authors focused on zero-inflated continuous distributions and proposed new estimators for Kendall's $\tau$ with reduced bias. 
\cite{denuit_bounds_2017} derived lower and upper bounds of the newly introduced estimator, making its interpretation possible as a measure of the strength of dependence.
The abundance of zero-inflated count data in practice, e.g., zero-inflated Poisson-type data, makes it crucial to define measures of dependence that can handle discreteness of the data as well as it being zero-inflated.
In this paper we extend the work of \cite{pimentel_association_2015} and propose a new estimator of Kendall's $\tau$ for bivariate random variables with zero-inflated discrete distributions. 
We complete the picture by deriving the theoretical lower and upper bounds of the proposed estimator, and we compare them with the bounds obtained by \cite{denuit_bounds_2017} for Pimentel's estimator. 
We show that our proposed estimator outperforms the available estimators in several simulated scenarios.
The paper is structured as follows: Section~\ref{sec:not} introduces the notation and basic concepts. In Section~\ref{sec:est}, we present our proposed estimator, and we discuss its attainable theoretical bounds in Section~\ref{sec:bounds}.
In Section~\ref{sec:sim}, we evaluate the performance of the estimator via a simulation study. We end with a discussion and conclusion section in Section~\ref{sec:concl}.

\section{Background and notation}
\label{sec:not}

We consider two independent copies $(\tilde{X}_1,\tilde{Y}_1)$ and $(\tilde{X}_2, \tilde{Y}_2)$ of the random vector $(X,Y)$ with joint cumulative distribution function $H$. Kendall's $\tau$ is defined as the probability of concordance minus the probability of discordance \citep{Kendall_1938}. 
For continuous random vectors, this definition results in 
$\tau= \mathbb{P}[(\tilde{X}_1 - \tilde{X}_2)(\tilde{Y}_1 - \tilde{Y}_2) > 0] - \mathbb{P}[(\tilde{X}_1 - \tilde{X}_2)(\tilde{Y}_1 - \tilde{Y}_2) < 0] = 2 \mathbb{P}[(\tilde{X}_1 - \tilde{X}_2)(\tilde{Y}_1 - \tilde{Y}_2) > 0] - 1.$
When $X$ and $Y$ assume values in the non-negative integers, Kendall's $\tau$ also depends on the probability of ties, i.e., $\tau=  \mathbb{P}[(\tilde{X}_1 - \tilde{X}_2)(\tilde{Y}_1 - \tilde{Y}_2) > 0] - 1 + \mathbb{P}[\tilde{X}_1 = \tilde{X}_2 \text{ or } \tilde{Y}_1=\tilde{Y}_2]$. 
The non-continuous case has extensively been studied in \cite{denuit_constraints_2002,Mesfioui_Tajar_2005, Neslehova_07}, and \cite{nikoloulopoulos_multivariate_2008}, where the authors also give closed-form formulas to calculate Kendall's $\tau$ for general discrete distributions when the distribution is completely known.

We denote as $\hat{\tau}$ the standard estimator of Kendall's $\tau$ computed by replacing the probability of concordance and discordance with the corresponding sample frequencies, which is by counting the number of concordant and discordant pairs and divide by the total number of pairs \citep{Kendall_1938}. 
In case of ties, i.e., repeated values in the sample, there are pairs that are neither concordant nor discordant. To account for this, an adjusted version of the estimator of Kendall's $\tau$ which excludes the tied pairs from the count has been proposed \citep{Kendall_1945}. We denote this by $\tau_b$. 
To define our theoretical framework, we use a similar notation as in \cite{pimentel_association_2015} and \cite{denuit_bounds_2017}.
We consider two non-negative random variables $X$ and $Y$ that follow two discrete distributions (e.g., Poisson) with extra positive probability mass at zero, i.e., the cumulative distribution function (cdf) of $X$, $F$, and of $Y$, $G$, can be written as follows
\begin{align*}
F(s) = \begin{cases} 0 \text{, if } s<0\\
(1 - \pi_F) + \pi_F \cdot \bar{F}(s) \text{, if } s \geq 0 
\end{cases} & \quad &
G(t) = \begin{cases} 0 \text{, if } t<0\\
(1 - \pi_G) + \pi_G \cdot \bar{G}(t) \text{, if } t \geq 0 
\end{cases}
\end{align*}
where, $\bar{F}$ and $\bar{G}$ are discrete distribution functions (e.g., Poisson), while for \cite{pimentel_association_2015} they are continuous.
We denote by $X_{10}$ a positive random variable distributed as $X$ given that $Y=0$, $X_{11}$ a positive random variable distributed as $X$ given that $Y>0$. Similarly, $Y_{01}$ is a positive random variable distributed as $Y$ given that $X=0$, and $Y_{11}$ a positive random variable distributed as $Y$ given that $X>0$.
We also consider $X_1$ a random variable distributed as $X$ given $X>0$ and $Y_1$, analogously.
In addition, we define the following probabilities: $p_{00} = \mathbb{P}[X=0, Y=0]$, $p_{01} = \mathbb{P}[X=0, Y>0]$, $p_{10} = \mathbb{P}[X>0, Y=0]$, $p_{11} = \mathbb{P}[X>0, Y>0]$, $p_1^{\ast} = \mathbb{P}[X_{10} > X_{11}]$, $p_2^{\ast} = \mathbb{P}[Y_{01} > Y_{11}]$, and define $\tau_{11}$ as Kendall's $\tau$ of $(X_1, Y_1)$, i.e., away from zero. 
Then, the association measure $\tau_H$ for the random vector $(X,Y)$ considered in \cite{pimentel_association_2015} is given by the following formula:
\begin{equation}
    \label{formula:Pim_Estimator}
    \tau_H = p_{11}^2 \tau_{11} + 2(p_{00}p_{11} - p_{01}p_{10}) + 2 p_{11}[p_{10}(1 - 2 p_1^{\ast}) + p_{01}(1- 2 p_2^{\ast})]
\end{equation}

\cite{pimentel_association_2015} suggested an estimator $\widehat{\tau_H}$ of Eq.~(\ref{formula:Pim_Estimator}) by replacing all probabilities of the formula with the corresponding sample frequencies, and, since no ties are expected away from zero, by substituting $\tau_{11}$ with the standard estimator $\hat{\tau}$ of Kendall's $\tau$ calculated from data on $X$ and $Y$ where $X$ and $Y$ are both positive.
In addition, \cite{denuit_bounds_2017} proved that the attainable bounds of the association measure $\tau_H$ only depend on the zero-inflation probabilities $p_1 = \mathbb{P}[X=0]$ and $p_2 =\mathbb{P}[Y=0]$. Specifically, the bounds are given by the following formulas
\begin{align}
\label{eq:pim_bounds_upper}
\tau^{upper}_H = & \begin{cases} 1 - p_2^2\text{, when } p_1 \leq p_2\\
1 - p_1^2\text{, when } p_1 \geq p_2
\end{cases}\\
\label{eq:pim_bounds_lower}
\tau^{lower}_H = & \begin{cases} - 2(1-p_1)(1-p_2)\text{, when } 1 - p_1 - p_2 <0
\\
(1 - p_1 - p_2)^2 - 2(1-p_1)(1-p_2)\text{, when } 1 - p_1 - p_2 >0
\\
\end{cases}
\end{align}
A natural question arises whether or not it is sufficient to replace the estimator of $\tau_{11}$ with $\tau_b$ in Eq.~(\ref{formula:Pim_Estimator}) to obtain an estimator of Kendall's $\tau$ to handle zero-inflated count data. The result established in the next section shows that this is not the case, and further adjustments are needed.

\section{Estimator of Kendall's $\tau$ for zero-inflated count data}
\label{sec:est}
The association measure studied in \cite{pimentel_association_2015} is based on a decomposition of the zero-inflated part from the continuous part of the distribution. 
\cite{pimentel_association_2015}'s estimator is interesting since, on the one hand, it accounts for the ties in zero and, on the other hand, it acts as the standard estimator of Kendall's $\tau$ away from zero.
Our approach is based on a similar idea of decomposing the association measure around zero and away from zero. 
However, due to the discrete nature of the zero-inflated count data away from zero, the estimator proposed by \cite{pimentel_association_2015} cannot be applied directly without further adjustments due to the non-zero probability of ties within the margins. 
The next result tackles this issue and establishes Kendall's $\tau$ for zero-inflated count data.

\begin{theorem}
\label{th:estimator}
We define the probabilities of ties within the margins as $p_1^{\dagger} = \mathbb{P}[X_{10}=X_{11}]$, and $p_2^{\dagger} = \mathbb{P}[Y_{01}=Y_{11}]$. Then Kendall's $\tau$ is given by the following relation
\begin{equation}
    \label{formula:est}
    \tau_A = p_{11}^2 \tau_{11} + 2(p_{00}p_{11} - p_{01}p_{10}) + 2 p_{11}[p_{10}(1 - 2 p_1^{\ast} - p_1^{\dagger}) + p_{01}(1- 2 p_2^{\ast} - p_2^{\dagger})].
\end{equation}
\end{theorem}
The proof of this theorem is straightforward and is similar to \cite{pimentel_kendalls_nodate}.
Based on the definition of Kendall's tau, we derived the expressions for the probability of concordance and discordance respectively using the law of total probabilities. The complete proof is given in the appendix.
As suggested in \cite{pimentel_association_2015}, an estimator $\widehat{\tau_A}$ of $\tau_A$ can be obtained by replacing probabilities with their estimates based on sample frequencies, while $\tau_{11}$ with the standard tie-corrected Kendall's $\tau$ estimator $\tau_b$ \citep{Kendall_1945}.
Moreover, consistency and asymptotic normality of the estimator $\widehat{\tau_A}$ follows directly from the same arguments presented in \cite{pimentel_association_2015} for the estimator of $\tau_H$.

\section{Attainable bounds for $\tau_A$}
\label{sec:bounds}
Kendall's $\tau$ cannot reach the theoretical bounds $\pm 1$ if there is a discrete component in the random vector.
In light of this, knowing the attainable bounds of the estimators of Kendall's $\tau$ is crucial to assess the strength of association of the data. 
To make our proposed estimator $\tau_A$ of Eq.~(\ref{formula:est}) useful in practice, we derive the range of admissible values for $\tau_A$ in terms of the marginal distributions of $X$ and $Y$. 
To do so, we follow the approach of \cite{denuit_bounds_2017}, which is based on the property of monotonicity of Kendall's $\tau$ with respect to the concordance order \citep{denuit_constraints_2002,Mesfioui_Tajar_2005}.

\begin{proposition}
\label{pr:bounds}
The lower and upper bounds of the association measure $\tau_A$ of Eq.~(\ref{formula:est}) are given by
\begin{align*}
\tau^{upper}_A = & \begin{cases} (1 - p_2^2) - (1 - p_2)^2p^U_{t_{11}} - 2(p_2 - F(\tilde{s}-1))(F(\tilde{s}) - p_2)  \text{, if } p_1 \leq p_2\\
(1 - p_1^2) - (1 - p_1)^2p^U_{t_{11}} - 2(p_1 - G(\tilde{t}-1))(G(\tilde{t}) - p_1) \text{, if } p_1 \geq p_2
\end{cases} \\[12pt]
\tau^{lower}_A = & \begin{cases} - 2(1-p_1)(1-p_2) \hfill \text{ if } 1 - p_1 - p_2 <0
\\[6pt]
p_1^2+p_2^2-1 + (1 - p_1 - p_2)^2 \cdot p^L_{t_{11}} + 2 [(F(\tilde{s}')+p_2-1)(1-p_2-F(\tilde{s}'-1)) + \\
\hspace{5.5cm} + (G(\tilde{t}')+p_1-1)(1-p_1-G(\tilde{t}'-1))], \hfill \text{ if } 1 - p_1 - p_2 >0 \\
\end{cases}
\end{align*}
where $\tilde{s}$ is a point such that $F(\tilde{s}) > p_2$ and $F(\tilde{s}-1) \leq p_2$, and $\tilde{s}'$ is a point such that $F(\tilde{s}') + p_2 -1 > 0$ and $F(\tilde{s}') + p_2 -1 \leq 0$ (analogously for $\tilde{t}$ and $\tilde{t}'$), and $p^U_{t_{11}}$ (or $p^L_{t_{11}}$) is the probability that either $X_1$ or $Y_1$ are tied when the joint distribution of $(X,Y)$ is the upper (lower) Fr\'echet-Hoeffding bound.
\end{proposition}

We notice that the points $\tilde{s}$, $\tilde{s}'$, $\tilde{t}$, $\tilde{t}'$ and the corresponding expressions are closely related to the joint probability expressed in terms of the Fr\'echet-Hoeffding bounds $\min\{F(x), G(y)\}$, and $\max\{F(x)+G(y)-1, 0\}$. 
The complete proof of Proposition~\ref{pr:bounds} is available in the appendix. 
Although the bounds reported in Proposition~\ref{pr:bounds} appear to be more involved than the bounds reported in \cite{denuit_bounds_2017}, it is still possible to estimate them from the data. The probabilities of zero-inflation of $X$ and $Y$, i.e., $p_1$ and $p_2$, can be replaced by the corresponding sample frequencies.
The values $p^U_{t_{11}}$ and $p^L_{t_{11}}$ are both dependent on the (joint) distribution. Nevertheless, as noticed in \cite{denuit_constraints_2002}, they can be replaced by their own lower bound $\max(\mathbb{P}[X_1 = X_2], \mathbb{P}[Y_1 = Y_2])$, which can be readily estimated from the sample without the knowledge of the distributions. 
Therefore, an estimator of a slightly wider range of $\tau_A$ can be constructed by substituting $p^U_{t_{11}}$ and $p^L_{t_{11}}$ with the maximum sample frequency of $X_1$ or $Y_1$ being tied.
Finally, the remaining values in the formulas, i.e., $F(\tilde{s}')$, $F(\tilde{s} -1)$, $F(\tilde{s}')$, $F(\tilde{s}'-1)$, and the analogous quantities for $G$, can be estimated via the empirical cdfs of $X$ and $Y$.

\section{Simulation study}
\label{sec:sim}
To investigate the performance of our proposed estimator, we conducted a Monte-Carlo simulation study based on 1000 repetitions. 
Namely, we computed the estimators of $\tau$, $\tau_H$ (with $\tau_{11}$ estimated via $\tau_b$ due to the discrete nature of our data), and $\tau_A$ for $N$ pairs generated from two correlated zero-inflated Poisson distributions joined through the Fr\'echet copula $C(u,v) = (1 - \rho) u v + \rho \min(u,v)$, where $u, v, \rho \in [0,1]$ \citep{nelsen_06}. 
Thus, the parameters of the full distribution are five, i.e., $\pi_F,\pi_G, \lambda_F, \lambda_G, \rho$, where $\rho$ is the copula parameter, $\pi_F$ is the probability mass spread according to a Poisson distribution with mean parameter $\lambda_F$ (analogously for $\pi_G$ and $\lambda_G$), and $1 - \pi_F$ represents the additional probability mass in zero which does not originate from the Poisson distribution of $X$ (same for $1-\pi_G$). 

We selected multiple scenarios representative of various characteristics of the samples.
In particular, we considered three values for the copula parameter $\rho=0.2, 0.5, 0.8$ depicting different strengths of association, two proportions $\pi_F = \pi_G = 0.2, 0.8$ for the probability mass associated with the Poisson distributions, and three combinations for the Poisson mean parameters $(\lambda_F, \lambda_G)=\{(2,2); (2,8) ; (8,8)\}$ corresponding to different levels of probability of ties away from zero.
We conducted our analysis in \texttt{R} \citep{RCore_2017} and made the code available as supplementary material.
We used the standard \texttt{R} function \texttt{cor()} to compute the tie-corrected version of $\tau$ (namely $\tau_b$), and we implemented the estimators of $\tau_A$ and $\tau_H$, which was also not available. 
In this regard, the meaning of the variable $X_{10}$ and $Y_{10}$ (respectively $X_{11}$ and $Y_{11}$) in \cite{pimentel_association_2015} are subject to interpretation.
In particular, in \cite{pimentel_association_2015}, $X_{10}$ is defined as a variable with a conditional distribution of $X$ given that $Y = 0$. 
Though, based on the steps of our proof and the notation chosen by the authors, we believe that $X_{10}$ should be a \emph{positive} random variable with a conditional distribution $X$ given that $Y = 0$ to result in a correct formulation for Eq.~(\ref{formula:Pim_Estimator}). 
To ensure a fair comparison with Pimentel's work, we implemented both interpretations for $X_{10}$ and selected the one that performs the best (which in fact corresponds to $X_{10}$ being a positive variable). 
We recall that, for a bivariate (zero-inflated) discrete distribution, we can calculate the true value of Kendall's $\tau$ as given in \cite{nikoloulopoulos_multivariate_2008}. Therefore, we compute the mean square error (MSE) of the estimators of $\tau$, $\tau_H$, and $\tau_A$ based on the true value of $\tau$, and identify the best one through the smallest MSE. Table~\ref{tab:mse_estimators} shows the results of our simulation study for our chosen parameter values and samples of size $N=150$.

\begin{figure}
\begin{centering}
\begin{tabular}{ccc}
  \includegraphics[scale=0.2]{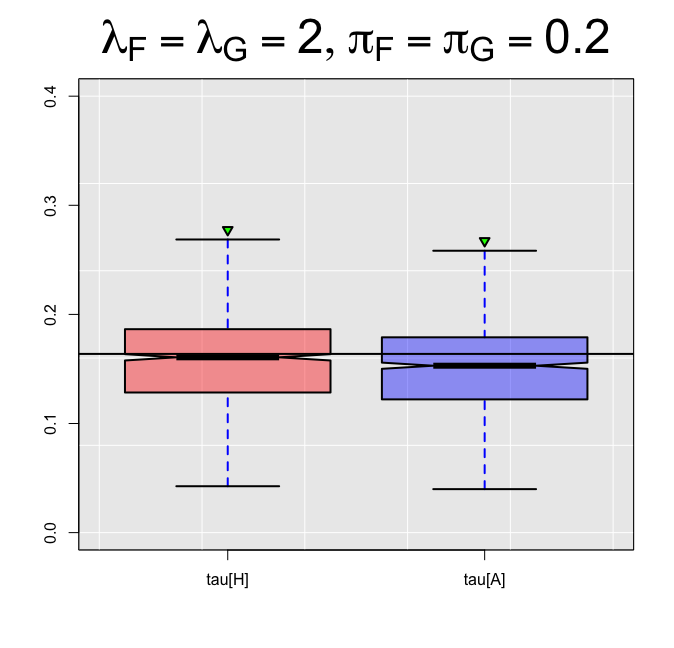} &  \includegraphics[scale=0.2]{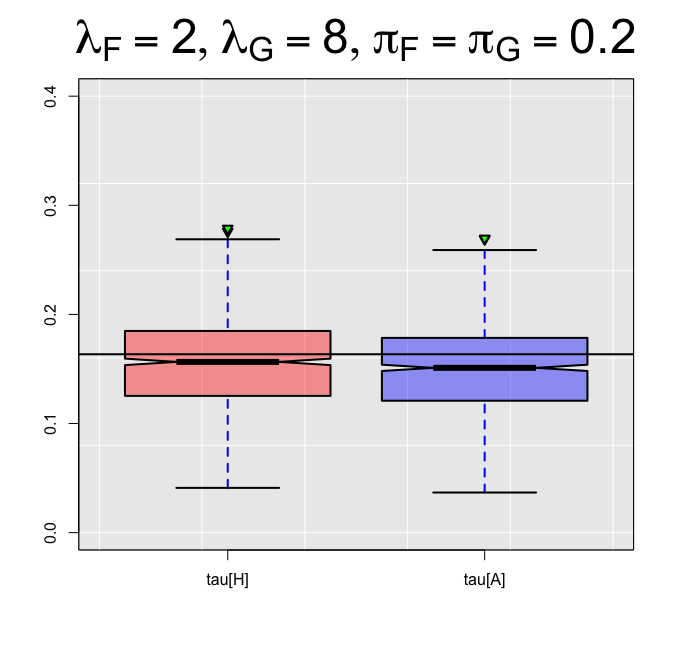} & \includegraphics[scale=0.2]{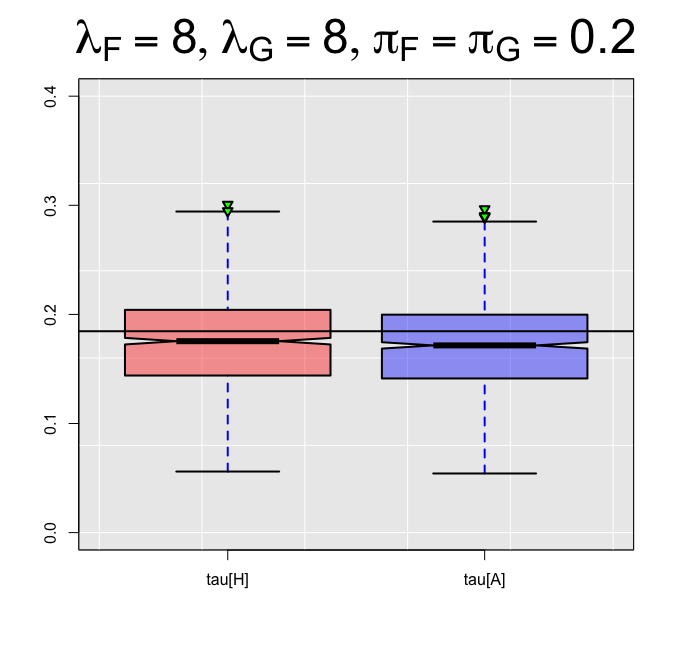} \\
  \includegraphics[scale=0.2]{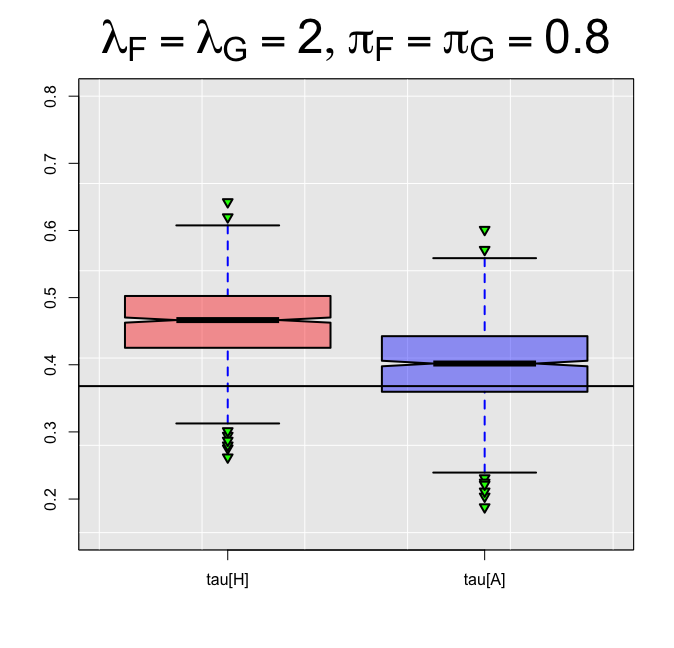} & 
  \includegraphics[scale=0.2]{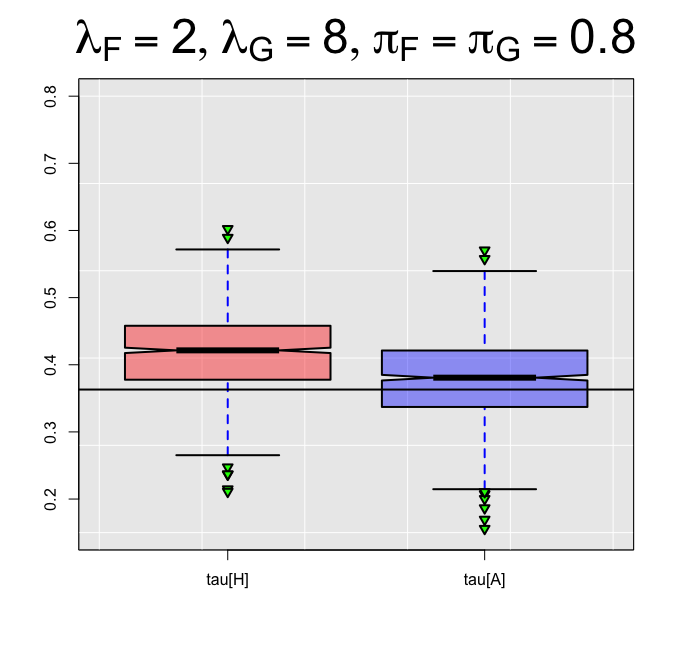}  &
   \includegraphics[scale=0.2]{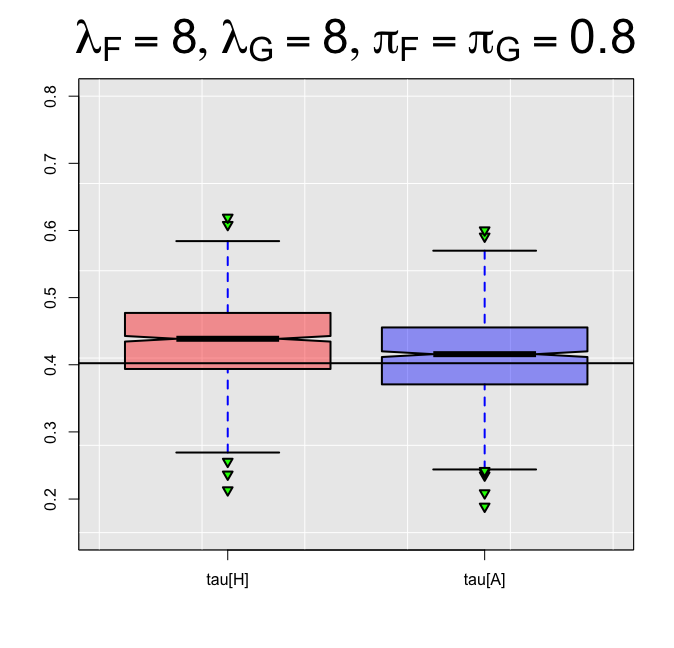}
    \\
\end{tabular}
\caption{Boxplots of $\tau_H$ (red) and $\tau_A$ (blue) over 1000 simulations for six different parameter settings and fixed $\rho=0.5$. The constant horizontal line in the plots represents the true value of $\tau$.}
\label{fig:boxplot_simulation}
\end{centering}
\end{figure}

\begin{table}[ht]
\centering
\begin{tabular}{|c|c|cc|lc|lc|}
\cline{5-8}
  \multicolumn{4}{c}{ } & \multicolumn{2}{|c|}{$\tau_H$ } & \multicolumn{2}{|c|}{$\tau_A$ } \\
  \hline
 & $\pi_F = \pi_G$ & $\rho$ & True $\tau$  &  Mean  & $\text{MSE}^*$ & Mean  & $\text{MSE}^*$ \\ 
\cline{2-8}
   &   & 0.20 & 0.07   & 0.07 & 0.11  & 0.06 & 0.12 \\ 
\multirow{3}{2.5cm}{$\lambda_F =2, \lambda_G=2$}  & 0.20  & 0.50 & 0.16  & 0.16 & 0.17  & 0.15 & 0.17 \\ 
  &  & 0.80 & 0.26  & 0.25 & 0.21  & 0.25 & 0.21 \\ 
\cline{2-8}
  &   & 0.20 & 0.15  & 0.24 & 1.16 & 0.15 & 0.38 \\ 
  & 0.80  & 0.50 & 0.37   & 0.46 & 1.23  & 0.40 & 0.48 \\ 
  &   & 0.80 & 0.62 &  0.72 & 1.21  & 0.69 & 0.77 \\ 
\hline
\hline
 &  & 0.20 & 0.07   & 0.06 & 0.12 & 0.06 & 0.12 \\ 
\multirow{4}{2.5cm}{$\lambda_F =2, \lambda_G=8$} & 0.20  & 0.50 & 0.16   & 0.16 & 0.17 & 0.15 & 0.18 \\ 
 &   & 0.80 & 0.26   & 0.25 & 0.21  & 0.25 & 0.21 \\ 
\cline{2-8}
 &   & 0.20 & 0.15  & 0.20 & 0.63  & 0.14 & 0.40 \\ 
 & 0.80  & 0.50 & 0.36   & 0.42 & 0.65  & 0.38 & 0.41 \\ 
 &   & 0.80 & 0.61   & 0.67 & 0.57  & 0.65 & 0.40 \\ 
   \hline
\hline
 &   & 0.20 & 0.08   & 0.07 & 0.13 &  0.07 & 0.14 \\ 
\multirow{4}{2.5cm}{$\lambda_F =8, \lambda_G=8$} & 0.20  & 0.50 & 0.18   & 0.18 & 0.19 & 0.17 & 0.19 \\ 
 &  & 0.80 & 0.29  & 0.28 & 0.23 &  0.28 & 0.23 \\ 
\cline{2-8}
 &  & 0.20 & 0.16   & 0.18 & 0.47 &  0.15 & 0.41 \\ 
 & 0.80  & 0.50 & 0.40  & 0.44 & 0.50 &  0.41 & 0.42 \\ 
 &  & 0.80 & 0.69  & 0.73 & 0.39 &  0.72 & 0.33 \\ 
   \hline
\end{tabular}
\caption{Comparison of the estimators' performance under various simulation scenarios and sample size $N=150$. The reported $\text{MSE}^*$ is the standard MSE multiplied by a factor of $10^2$.}
\label{tab:mse_estimators}
\end{table}
We investigated various sample sizes, i.e., $N=150, 300, 1000$, and did not find differences in the behavior of the estimators. Thus, we only present the case $N=150$ in the paper. Given the poor performances of $\tau_b$ in all the considered scenarios, we decided not to report it in the table.
As expected, the estimators of $\tau_H$ and $\tau_A$ are comparable in their performance when (1) most of the probability mass is in zero for both variables, i.e., for $\pi_F = \pi_G = 0.2$, and (2) there is a limited proportion of ties within the margins away from zero, i.e., $(\lambda_F, \lambda_G) = (8,8)$. 
In the other cases, our proposed estimator $\tau_A$ outperforms $\tau_H$. 
A visual representation of the performance of the estimators for selected parameter settings is presented in Figure~\ref{fig:boxplot_simulation}. From the boxplots of Figure~\ref{fig:boxplot_simulation}, we can conclude that $\tau_A$ is generally close to the true value of $\tau$ (the constant horizontal line in the plots), while $\tau_H$ tends to overestimate it when the zero-inflation is mild (e.g., $\pi_F=\pi_G=0.8$). 
Collectively, our analysis demonstrates that replacing $\tau_{11}$ in Eq.~(\ref{formula:Pim_Estimator}) by $\tau_b$ is not enough to ensure accurate performances, and our adjustment accounting for the probability of ties within the margins is needed when dealing with zero-inflated count data. 

Besides looking at the performance of the estimators, we also investigated their attainable bounds. We estimated the range of $\tau_H$ as suggested in \cite{denuit_bounds_2017} and the bounds of $\tau_A$ as described in Section~\ref{sec:bounds}. 
For comparison, we applied Proposition~\ref{pr:bounds} and find the theoretical bounds of $\tau_A$ for our specific marginal distributions (Zero-Inflated Poisson).
The results are reported in Table~\ref{tab:bounds_estimators}. 
The estimated ranges for $\tau_H$ and $\tau_A$ are very similar to each other and close to the theoretical bounds of $\tau_A$ when a small probability mass is spread away from zero. 
When zero-inflation is limited, our estimated bounds are sharper than the ones derived in \cite{denuit_bounds_2017}, as expected. 
We notice that they are still loose if compared to the true range as they are not exact estimates of the bounds of Proposition~\ref{pr:bounds}.
Nevertheless, such a non-parametric estimator of the bounds is very useful in practice as it allows for interpreting the strength of association of an estimate of $\tau_A$ without the need to make assumptions on the underlying distributions.

\begin{table}[h]
\centering
\begin{tabular}{|c|c|ccc|}
  \hline
\multirow{3}{2.5cm}{$\lambda_F =2, \lambda_G=2$} & $\pi_F = \pi_G$ &  Bounds $\tau_H$  & Bounds $\tau_A$ & Theoretical bounds $\tau_A$\\
\cline{2-5}
 &  0.20   & [-0.06, 0.29]  & [-0.06, 0.29] & [-0.06, 0.31] \\ 
 & 0.80   & [-0.81, 0.90]  & [-0.76, 0.84] & [-0.75, 0.78] \\ 
   \hline
   \hline
\multirow{2}{2.5cm}{$\lambda_F =2, \lambda_G=8$}   &  0.20  & [-0.07, 0.32]  & [-0.07, 0.32] & [-0.07, 0.31] \\ 
  & 0.80   & [-0.86, 0.90]  & [-0.82, 0.85] & [-0.80, 0.77] \\ 
   \hline
   \hline
\multirow{2}{2.5cm}{$\lambda_F =8, \lambda_G=8$} &     0.20   & [-0.06, 0.29]  & [-0.06, 0.29] & [-0.06, 0.31] \\ 
 & 0.80    & [-0.81, 0.90] & [-0.76, 0.84] & [-0.75, 0.78] \\
 \hline
\end{tabular}
\caption{Estimates of the lower and upper bounds of $\tau_H$ and $\tau_A$ for various simulation scenarios, sample size $N=150$, and averaged across 1000 runs.}
\label{tab:bounds_estimators}
\end{table}

\section{Conclusion}
\label{sec:concl}
In this paper, we built on previous results presented in \cite{pimentel_association_2015} by proposing an adjusted estimator of Kendall's $\tau$ that can tackle both zero-inflated continuous and count data.
We made the proposed estimator interpretable and useful in practice by deriving its theoretical attainable bounds and suggesting a way to estimate them. 
Our theoretical results were paired with a simulation study, where we analyzed the estimators' performance in various settings. 
Overall, our proposed estimator is more flexible and preferable in practice since it coincides with the estimator proposed by \cite{pimentel_association_2015} if there are no ties within the margins, while it outperforms it if the zero-inflation is mild.

This paper was motivated by the need for \emph{ad hoc} statistical methods to quantify association between zero-inflated count data. 
Thus, a natural follow-up of this work would be to derive a suitable estimator for Spearman's $\rho$ in case of zero-inflated count data. 
Preliminary results on the topic for the continuous zero-inflated case have been presented in \cite{pimentel_kendalls_nodate, Mesfioui2022a}, while \cite{Mesfioui2022c} recently analyzed the attainable bounds of Spearman's $\rho$ when at least one variable is discrete. 
However, more research and investigations are needed to derive for Spearman's $\rho$ the same tools now available for Kendall's $\tau$. 

\bibliography{CEBZIP.bib}

\section*{Appendices}

\section*{A1. Proof of Theorem~\ref{th:estimator}}

For the two independent and identically distributed sample $(\xx_1, \yy_1)$ and $(\xx_2, \yy_2)$ of the random vector $(X, Y)$ with joint cumulative distribution function $H$, the definition of Kendall's $\tau$ can be written as
\begin{equation*}
  \tau_A = \pp(\xx_1 < \xx_2, \yy_1 < \yy_2) + \pp(\xx_1 > \xx_2, \yy_1 > \yy_2) - \pp(\xx_1 < \xx_2, \yy_1 > \yy_2) - \pp(\xx_1>\xx_2, \yy_1 < \yy_2).
\end{equation*}
Furthermore, consider the following set $C_{a_{1}b_{1}a_{2}b_{2}} = \{\sign \xx_1 = a_1, \sign \yy_1 = b_1, \sign \xx_2 =a_2, \sign \yy_2 =b_2\}$, with $a_1, b_1, a_2, b_2 \in \{0,1\}$, and $\sign x$ the sign function. Then, for instance, $C_{1111}$ denotes the set $\{\xx_1 >0 ,\yy_1 >0, \xx_2 >0, \yy_2 >0\}$ and $\pp(C_{1111})=p_{11}p_{11}$. More generally, $\pp(C_{a_{1}b_{1}a_{2}b_{2}}) = p_{a_1 b_1}p_{a_2 b_2}$. Using the law of total probability, it can be shown that
\begin{equation*}
\begin{aligned}
\pp(\xx_1 < \xx_2, \yy_1 < \yy_2) = & \sum_{(a_1, b_1, a_2, b_2) \in \{0,1\}^4} \pp(\xx_1 < \xx_2, \yy_1 < \yy_2 | C_{a_{1}b_{1}a_{2}b_{2}}) \pp(C_{a_{1}b_{1}a_{2}b_{2}})\\
= & \pp(X_1 < X_2, Y_1 < Y_2 | C_{1111}) p_{11}^2 + \pp(X_1 < X_2, Y_1 < Y_2 | C_{1011}) p_{10}p_{11} + \\
 & \pp(X_1 < X_2, Y_1 < Y_2 | C_{0111}) p_{01}p_{11} + p_{00}p_{11}.\\
\end{aligned}
\end{equation*}
Similarly, we can rewrite the other three probabilities in the right hand side of the definition. This yields the following expressions for the four terms, respectively:
\begin{equation*}
  \begin{aligned}
  &\pp(\xx_1 < \xx_2, \yy_1 < \yy_2)\\
  =&\pp(\xx_1 < \xx_2, \yy_1 < \yy_2 | C_{1111}) p_{11}^2 +\pp(\xx_1 < \xx_2, \yy_1 < \yy_2 | C_{1011}) p_{10}p_{11} + \pp(\xx_1 < \xx_2, \yy_1 < \yy_2 | C_{0111}) p_{01}p_{11} + p_{00}p_{11},\\[6pt]
   &\pp(\xx_1 > \xx_2, \yy_1 > \yy_2) \\
  =&\pp(\xx_1 > \xx_2, \yy_1 > \yy_2 | C_{1111}) p_{11}^2 + \pp(\xx_1 > \xx_2, \yy_1 > \yy_2 | C_{1101}) p_{11}p_{01} + \pp(\xx_1 > \xx_2, \yy_1 > \yy_2 | C_{1110}) p_{11}p_{10} + p_{00}p_{11},\\[6pt]
   &\pp(\xx_1 < \xx_2, \yy_1 > \yy_2) \\
  =&\pp(\xx_1 < \xx_2, \yy_1 > \yy_2 | C_{1111}) p_{11}^2 + \pp(\xx_1 < \xx_2, \yy_1 > \yy_2 | C_{0111}) p_{01}p_{11} + \pp(\xx_1 < \xx_2, \yy_1 > \yy_2 | C_{1110}) p_{11}p_{10} + p_{01}p_{10},\\[6pt]
  &\pp(\xx_1 > \xx_2, \yy_1 < \yy_2) \\
 =&\pp(\xx_1 > \xx_2, \yy_1 < \yy_2 | C_{1111}) p_{11}^2 + \pp(\xx_1 > \xx_2, \yy_1 < \yy_2 | C_{1101}) p_{11}p_{01} + \pp(\xx_1 > \xx_2, \yy_1 < \yy_2 | C_{1011}) p_{10}p_{11} + p_{10}p_{01}.\\[6pt]
  \end{aligned}
\end{equation*}
Consequently, the Kendall's $\tau$ can be expressed as
\begin{equation*}
  \tau_A = p_{11}^2 \tau_{11} + p_{11} A + 2(p_{00}p_{11} - p_{01}p_{10})
\end{equation*}
where 
\begin{equation*}
  \tau_{11} =\pp(\xx_1 < \xx_2, \yy_1 < \yy_2 | C_{1111})  + \pp(\xx_1 > \xx_2, \yy_1 > \yy_2 | C_{1111}) - \pp(\xx_1 < \xx_2, \yy_1 > \yy_2 | C_{1111})  - \pp(\xx_1 > \xx_2, \yy_1 < \yy_2 | C_{1111})
\end{equation*}
is the Kendall's $\tau$ conditioned on $C_{1111}$, and
\begin{equation*}
  \begin{aligned}
    A &=  \pp(\xx_1 < \xx_2, \yy_1 < \yy_2 | C_{1011}) p_{10} + \pp(\xx_1 < \xx_2, \yy_1 < \yy_2 | C_{0111}) p_{01} \\
    &+\pp(\xx_1 > \xx_2, \yy_1 > \yy_2 | C_{1101}) p_{01} + \pp(\xx_1 > \xx_2, \yy_1 > \yy_2 | C_{1110}) p_{10}\\
    &-\pp(\xx_1 < \xx_2, \yy_1 > \yy_2 | C_{0111}) p_{01} - \pp(\xx_1 < \xx_2, \yy_1 > \yy_2 | C_{1110}) p_{10} \\
    &-\pp(\xx_1 > \xx_2, \yy_1 < \yy_2 | C_{1101}) p_{01} -  \pp(\xx_1 > \xx_2, \yy_1 < \yy_2 | C_{1011}) p_{10}.
  \end{aligned}
\end{equation*}
Recall that $(\xx_1, \yy_1)$ and $(\xx_2, \yy_2)$ are exchangeable by definition, therefore the terms in the expression of $A$ can be further simplified as
\begin{equation*}
\begin{aligned}
  A = &2\Bigl\{p_{10}\bigl[\pp(\xx_1 < \xx_2, \yy_1 < \yy_2 | C_{1011}) - \pp(\xx_1 > \xx_2, \yy_1 < \yy_2|C_{1011})\bigr] \\
  &+ p_{01}\bigl[\pp(\xx_1 < \xx_2, \yy_1 < \yy_2|C_{0111}) - \pp(\xx_1 < \xx_2, \yy_1 > \yy_2|C_{0111})\bigr]\Bigr\}\\
  =& 2\Bigl\{p_{10}\bigl[\pp(\xx_1 < \xx_2 | C_{1011}) - \pp(\xx_1 > \xx_2 |C_{1011})\bigr] + p_{01}\bigl[\pp(\yy_1 < \yy_2|C_{0111}) - \pp(\yy_1 > \yy_2 |C_{0111})\bigr]\Bigr\} \\
 =&2\Bigl\{p_{10}\bigl[1 - 2\pp(\xx_1 > \xx_2 | C_{1011}) - \pp(\xx_1 = \xx_2 |C_{1011})\bigr] + p_{01}\bigl[1 - 2\pp(\yy_1 > \yy_2|C_{0111}) - \pp(\yy_1 = \yy_2 |C_{0111})\bigr]\Bigr\}\\
 =&2\bigl[p_{10}(1-2p_1^\ast - p_1^\dagger) + p_{01}(1 - 2p_2^\ast - p_1^\dagger)\bigr],
\end{aligned}
\end{equation*}
with $p_i^\ast$ and $p_i^\dagger$ for $i=1,2$ as defined in the main text. Finally, we obtain the expression for Kendall's $\tau$ given by
\begin{equation*}
  \tau_A = p_{11}^2 \tau_{11} + 2(p_{00}p_{11} - p_{01}p_{10}) + 2p_{11}[p_{10}(1-2p_1^\ast - p_1^\dagger) + p_{01}(1 - 2p_2^\ast - p_1^\dagger)].
\end{equation*}
\qed

\section*{A2. Distributions of the random variables $X_{10}$ and $X_{11}$ under the Fr\'echet-Hoeffding bounds}
In this section, we try to establish the distributions of the random variables $X_{10}$ and $X_{11}$ defined in the main text when the joint distribution $H$ of the random pair $(X, Y)$ is equal to the Fr\'echet-Hoeffding bounds.

When $H(x,y) = \min\{F(x), G(y)\}$, we can write the cumulative probability function of $X_{10}$ as
\begin{equation*}
  \begin{aligned}
    \pp(X_{10} \le x) &= \pp(X \le x | X >0, Y = 0) \\
    &=\frac{\pp(0 < X \le x, Y=0)}{\pp(X>0, Y = 0)}\\
    &= \frac{1}{p_2 - p_1} \sum_{s=1}^x \pp(X=s, Y=0)\\
    &= \frac{1}{p_2 - p_1}\sum_{s=1}^x \bigl\{\pp(X \le s, Y \le 0) - \pp(X\le s-1, Y\le 0)\bigr\}\\
    &=\frac{1}{p_2 - p_1} \sum_{s=1}^x \min\{F(s), p_2\} - \min\{F(s-1), p_2\}.
  \end{aligned}
\end{equation*}
Let $\ss$ denotes the point in the support of $F(x)$ such that $F(\ss) > p_2$ and $F(\ss-1) \le p_2$, the probability mass function of $X_{10}$ is as follows
\begin{equation*}
  P(X_{10} = x) = \begin{cases} 
    \frac{F(x) - F(x-1)}{p_2 - p_1} & 0 < x < \ss \\
    \frac{p_2 - F(x-1)}{p_2 - p_1} & x = \ss \\
    0 & x > \ss
  \end{cases}
\end{equation*}
In a similarly manner, it can be shown for the random variable $X_{11}$ that its cumulative distribution and probability mass function is given by
\begin{equation*}
  \begin{aligned}
    \pp(X_{11} \le x) &= \pp(X \le x | X >0, Y>0) \\
    & = \frac{1}{1-p_2}\bigl\{\pp(0 < X \le x) - \pp(0 < X \le x, Y=0)\bigr\} \\
    & = \frac{1}{1-p_2} \sum_{s=1}^x \pp(X=s) - \pp(X=s, Y=0),
  \end{aligned}
\end{equation*}
and, respectively,  
\begin{equation*}
  P(X_{11} = x) = \begin{cases} 
    0 & 0 < x < \ss \\
    \frac{F(x) - p_2}{1-p_2} & x = \ss \\
    \frac{F(x) - F(x-1)}{1 - p_2} & x > \ss
  \end{cases}.
\end{equation*}
Furthermore, we can calculate the probability $\pp(X_{10} > X_{11})$ as
\begin{equation*}
  \pp(X_{10} > X_{11}) = \sum_{x=0}^\infty \pp(X_{10} > x)\pp(X_{11} = x) = \sum_0^\infty \bigl\{1 - \pp(X_{10} \le x)\bigr\} \pp(X_{11} = x).
\end{equation*}
It can be seen that when $x < \ss$, the second term of the product in the summation is 0, and when $x \ge \ss$, the first term becomes 0. Therefore, $\pp(X_{10} > X_{11}) = 0$. On the other hand, for 
\begin{equation*}
  \pp(X_{10} = X_{11}) = \sum_{x=0}^\infty \pp(X_{10}=x)\pp(X_{11}=x) = \frac{p_2 - F(\ss-1)}{p_2 - p_1}\frac{F(\ss) - p_2}{1 - p_2}
\end{equation*}
since the summation has all its summand equal to 0 except when $x=\ss$.

When $H(x,y) = \max\{F(x) + G(y) - 1, 0\}$, we only consider the case when $p_1 + p_2 -1 \le 0$. The cumulative distribution function of $X_{10}$ in this situation can be written as
\begin{equation*}
  \begin{aligned}
    \pp(X_{10} \le x) & = \pp(X \le x | X >0, Y=0) \\
    &=\frac{\pp(0 < X \le x, Y=0)}{\pp(X >0, Y=0)}\\
    &=\frac{1}{p_2} \sum_{s=1}^x \pp(X=s, Y=0) \\
    &=\frac{1}{p_2} \sum_{s=1}^x\bigl\{\pp(X \le s, Y=0) - \pp(X\le s-1, Y=0)\bigr\},
  \end{aligned}
\end{equation*}
and the probability mass function is 
\begin{equation*}
  \pp(X_{10} = x) = \begin{cases}
    0 & x < \ss' \\
    \frac{F(\ss') + p_2 - 1}{p_2} & x = \ss' \\
    \frac{F(x) - F(x-1)}{p_2} & x > \ss'
  \end{cases}
\end{equation*}
where $\ss'$ is a point in the support of $F(x)$ such that $F(\ss') + p_2 - 1 > 0$ and $F(\ss' - 1) + p_2 - 1 \le 0$.
Similarly, $X_{11}$ has a cumulative distribution function of
\begin{equation*}
  \pp(X_{11} \le x) = \frac{1}{1 - p_1 - p_2} \sum_{s=1}^x \bigl\{\pp(X=s) - \pp(X=s, Y=0)\bigr\},
\end{equation*}
and a probability mass function that is given by
\begin{equation*}
  \pp(X_{11} = x) = \begin{cases}
    \frac{F(x) - F(x-1)}{1 - p_1 - p_2} & x < \ss' \\
    \frac{1 - p_2 - F(\ss' - 1)}{1 - p_1 - p_2} & x = \ss' \\
    0 & x > \ss'
  \end{cases}.
\end{equation*}
Accordingly, the probabilities $\pp(X_{10} \le X_{11})$ and $\pp(X_{10} = X_{11})$ are given by
\begin{equation*}
  \pp(X_{10} \le X_{11}) = \sum_{x=0}^\infty \pp(X_{10} \le x)\pp(X_{11} = x) = \frac{F(\ss') + p_2 - 1}{p_2}\frac{1-p_2-F(\ss'-1)}{1-p_1-p_2}
\end{equation*}
and, respectively,
\begin{equation*}
  \pp(X_{10} = X_{11}) = \sum_{x=0}^\infty \pp(X_{10} = x)\pp(X_{11} = x) = \frac{F(\ss') + p_2 - 1}{p_2}\frac{1-p_2-F(\ss'-1)}{1-p_1-p_2}.
\end{equation*}
This shows that $\pp(X_{10} \le X_{11}) = \pp(X_{10} = X_{11})$ and $\pp(X_{10} < X_{11})=0$.

\section*{A3. Proof of Proposition~\ref{pr:bounds}}
Following \cite{denuit_bounds_2017}, the upper and lower bound of $\tau_A$ can be derived based on the monotonicity property of Kendall's $\tau$. Namely, for two joint distribution functions $H_1$ and $H_2$, if $H_1(x, y) \le H_2(x, y)$ for all $x$ and $y$, then $\tau_{H_1} \le \tau_{H_2}$.

Our proposed $\tau_A$ attains its upper bound $\tau_{A}^{upper}$ when $H(x,y) = \min\{F(x), G(y)\}$. Without loss of generality, we assume $p_1 \le p_2$. Consequently, $p_{00} = p_1, p_{01} = 0, p_{10} = p_2 - p_1$, and $p_{11} = 1- p_2$. The upper bound $\tau_{A}^{upper}$ is given by
\begin{equation*}
  \tau_{A}^{upper} = p_{11}^2 \tau_{11}^U + 2p_{00}p_{11} + 2p_{11}p_{10}\bigl\{1 - \pp(X_{10} = X_{11})\bigr\}
\end{equation*}
where $\tau_{11}^U$ is an upper bound of the $\tau_{11}$. Here, we take $\tau_{11}^U = 1 - p_{t_{11}}^U$ with $p_{t_{11}}^U$ the probability that either $X_1$ or $Y_1$ are tied when $H(x,y) = \min\{F(x), G(y)\}$. Based on the distributions of $X_{10}$ and $X_{11}$ in the previous section, the upper bound equals to
\begin{equation*}
  \tau_A^{upper} = (1-p_2^2) - (1-p_2)^2p_{t_{11}}^U - 2(p_2 - F(\ss -1))(F(\ss) - p_2).
\end{equation*}

When $H(x,y) = \max\{F(x) + G(y) -1, 0\}$ and $1 - p_1 - p_2 \ge 0$, we have $p_{00} = 0$, $p_{01} = p_1$, $p_{10} = p_2$, $p_{11} = 1 - p_1 - p_2$, and the lower bound $\tau_A^{lower}$ of $\tau_A$ is given by
\begin{equation*}
  \begin{aligned}
    \tau_A^{lower} &= p_{11}^2 \tau_{11}^L + 2(p_{00}p_{11} - p_{01}p_{10}) + 2p_{11}\Bigl\{-p_{10}\pp(X_{10} > X_{11}) - p_{01}\pp(Y_{10} > Y_{11})\Bigr\} \\
    &= p_1^2 + p_2^2 -1 + (1 - p_1 - p_2)^2 p_{t_{11}}^L + 2[(F(\ss') + p_2 - 1)(1 - p_2 - F(\ss' - 1)) +\\
    & \quad (G(\tt') + p_1 -1)(1-p_1-G(\tt' - 1))].
  \end{aligned}
\end{equation*}
The second equality follows when substituting the $\tau_{11}$'s lower bound $\tau_{11}^L$ with $p_{t_{11}}^L - 1$ where $p_{t_{11}}^L$ is defined analogously as $p_{t_{11}}^U$.

When $H(x,y) = \max\{F(x) + G(y) -1, 0\}$ and $1 - p_1 - p_2 < 0$, then $p_{00} = p_1 + p_2 - 1$, $p_{01} = 1 - p_2$, $p_{10} = 1 - p_1$, and $p_{11} = 0$. Therefore, the lower bound $\tau_A^{lower} = -2(1-p_1)(1-p_2)$.

\qed

\end{document}